\documentclass[12pt,a4paper]{article}

\usepackage[latin2]{inputenc}
\usepackage{amssymb}
\usepackage{paralist}
\usepackage{amsmath}
\usepackage{amsfonts}
\usepackage{amsthm}
\usepackage{fullpage}
\usepackage{enumerate}
\newtheorem{thm}{Theorem}
\newtheorem{lem}{Lemma}
\newtheorem{prob}{Problem}
\newtheorem{conj}{Conjecture}
\newtheorem{cor}{Corollary}

\begin{document}

\title{On the maximum values of the additive representation functions}
\author{S\'andor Z. Kiss \thanks{Institute of Mathematics, Budapest
University of Technology and Economics, H-1529 B.O. Box, Hungary;
Computer and Automation Research Institute of the Hungarian Academy of
Sciences, Budapest H-1111, L\'agym\'anyosi street 11;
kisspest@cs.elte.hu; MTA-BME Lend\"ulet Future Internet Research Group;
This author was supported by the OTKA Grant No. NK105645.}, Csaba
S\'andor \thanks{Institute of Mathematics, Budapest University of
Technology and Economics, H-1529 B.O. Box, Hungary, csandor@math.bme.hu.
This author was supported by the OTKA Grant No. K109789. This paper was supported
by the J\'anos Bolyai Research Scholarship of the Hungarian Academy of Sciences.}
}
\date{}
\maketitle

\begin{abstract}
\noindent Let $A$ and $B$ be sets of nonnegative integers. For a positive
integer $n$ let $R_{A}(n)$ denote the number of representations of $n$ as the
sum of two terms from $A$. Let $\displaystyle s_{A}(x) = \max_{n \le
 x}R_{A}(n)$ and $\displaystyle d_{A,B}(x) = \max_{\hbox{t: $a_{t} \le x$ or $b_{t} \le x$}}|a_{t} - b_{t}|$. In this paper we study the connection between
 $s_{A}(x)$, $s_{B}(x)$ and $d_{A,B}(x)$. We improve a result of Haddad
 and Helou about the Erd\H{o}s - Tur\'an conjecture.
\end{abstract}
\textit{2000 AMS \ Mathematics subject classification number}: 11B34.

\section{Introduction}
Let $\mathbb{N}$ denote the set of nonnegative integers. Let $A =
\{a_{1}, a_{2},\ldots{}\}$, $0 \le a_{1} < a_{2} < \dots{}$, $B =
\{b_{1}, b_{2},\ldots{}\}$, $0 \le b_{1} < b_{2} < \dots{}$ be infinite
sequences of nonnegative integers.  Let $R_{A}(n)$
denote the number of solutions of the equation
\[
a_{i} + a_{j} = n, \hspace*{3mm} a_{i}, a_{j} \in A,
\hspace*{3mm} i \le j,
\]
\noindent where $n \in \mathbb{N}$. Let $A(x)$ denote the number of elements of the sequence $A$ up to $x$. A (finite or infinite) set $A$ of
positive integers is said to be a Sidon set if all the sums $a + b$ with
$a,b \in A$, $a \le b$  are distinct.
Let $\displaystyle s_{A}(x) = \max_{n \le x}R_{A}(n)$ and $\displaystyle s_{A} = \lim_{x \rightarrow
  \infty}s_{A}(x)$, $\displaystyle d_{A,B}(x) =\max_{\hbox{t: $a_{t} \le x$ or $b_{t} \le x$}}|a_{t} - b_{t}|$, $\displaystyle d_{A,B} = \lim_{x \rightarrow \infty}d_{A,B}(x)$.
Throughout the paper we use the following usual notations. If $f$ and $g$ are real functions, then $f \ll g$ means that $f = O(g)$. If there exist positive constants $c_{1}$ and $c_{2}$ such that $c_{1}g \le f \le c_{2}g$ then we write $f = \Theta(g)$. We write $f \sim g$ if the function $f$ is asymptotically equal to $g$.

In 1941 Erd\H{o}s and Tur\'an [7] posed the following conjecture: if
$R_{A}(n)$ is positive
from a certain point on, then it cannot be bounded. Despite all the efforts
this conjecture is still unsolved.
In [8] and [9] Erd\H{o}s and Fuchs formulated the following conjecture
which is a little bit stronger then the original conjecture of Erd\H{o}s
and Tur\'an.
\begin{conj}
For any $\mathcal{A} = \{a_{1}, a_{2},\ldots{}\}$ set of nonnegative integers
with $a_{n} \le cn^{2}$ for all $n$ and $c > 0$ real constant, we have
 $\displaystyle \limsup _{n\to \infty}R_{A}(n) = \infty$.
\end{conj}
\noindent It is clear that this conjecture implies the original conjecture of
Erd\H{o}s and Tur\'an. However, a lot of partial results has been made
about the Erd\H{o}s - Tur\'an conjecture [2, 3, 4, 5, 6, 7, 8, 9, 10,
15, 16, 18, 19, 20, 21, 22, 24, 25, 26, 28],
very little progress has been made about the generalized Erd\H{o}s -
Tur\'an conjecture. In [10], [12] Grekos,
Haddad, Helou and Pihko proved several statements that are equivalent to the
generalized Erd\H{o}s - Tur\'an conjecture. In particular, they proved
[12], [13] that
the conjecture is true if $a_{n} = o(n^{2})$. In [16] Haddad and Helou
proved the following theorem:

\begin{thm}[Haddad, Helou, 2012]
If $|a_{n} - n^{2}| = o(\sqrt{\log n})$, or in general
$|a_{n} - q(n)| = o(\sqrt{\log n})$, where $q(n)$ is a quadratic polynomial
with rational coefficients, then $R_{A}(n)$ is unbounded.
\end{thm}
\noindent In this paper we improve on their result by proving that
\begin{thm}
For an arbitrary $\varepsilon > 0$ if $A \subset \mathbb{N}$, $A = \{a_{1}, a_{2}, \dots{}\}$ such that
 $|a_{n} - n^{2}| = O\Big(e^{(\log 2 - \varepsilon)\frac{\log n}{\log \log
     n}}\Big)$, then $R_{A}(n)$ cannot be bounded.
\end{thm}
\noindent We conjecture that this result can be generalized.
\begin{prob}
Let $q(n)$ be a real quadratic polynomial with positive leading coefficient. Is it true that for $|a_{n} - q(n)| = O\Big(e^{(\log 2 - \varepsilon)\frac{\log n}{\log \log n}}\Big)$ we have $s_A=\infty$?
\end{prob}
\noindent The aim of this paper is to study the connections between the
functions $s_{A}(x)$, $d_{A,B}(x)$ and $s_{B}(x)$. We distinguish four cases
according to $s_{A}$ and $d_{A,B}$ are finite or infinite. Theorem 3. deals
with the case $s_{A}$ and $d_{A,B}$ are finite. We determine all the possible
triplets $(s_{A}$, $s_{B}$, $d_{A,B})$. Note that the first part of Theorem 3.
is Corollary 3. in [12].

\begin{thm}
\begin{enumerate}
\item Let $A, B \subset \mathbb{N}$ such that $s_{A} < \infty$, $d_{A,B} <
\infty$. Then we have
\[
\frac{s_{A}}{4d_{A,B}+1} \le s_{B} \le (4d_{A,B}+1)s_{A}.
\]
\item Let $a, b, d \in \mathbb{Z^{+}}$ such that
\[
\frac{a}{4d+1} \le b \le (4d+1)a.
\]
Then there exist $A, B \subset \mathbb{N}$ such that $s_{A} = a$, $s_{B} = b$,
$d_{A,B} = d$.
\end{enumerate}
\end{thm}
\noindent The following theorem deals with the case when the representation
function $R_A(n)$ is unbounded, but the distance $d_{A,B}(n)$ is bounded.
\begin{thm}
\begin{enumerate}
\item Let $A$, $B \subset \mathbb{N}$ such that $s_{A} = +\infty$, $d_{A,B} <
+\infty$. Then
\[
\frac{1}{4d_{A,B} + 1} \le \liminf_{x \rightarrow \infty}
\frac{s_{B}(x)}{s_{A}(x - 2d_{A,B})} \le  \limsup_{x \rightarrow \infty}
\frac{s_{B}(x)}{s_{A}(x + 2d_{A,B})} \le 4d_{A,B} + 1.
\]
\item Let $d$ be a positve integer and $\alpha$, $\beta$ positive real numbers such
that $\frac{1}{4d + 1} \le \alpha \le \beta \le 4d + 1$. Then there exist $A$,
$B \subset \mathbb{N}$ such that $s_{A} = +\infty$, $d_{A,B} = d$ and
\[
\liminf_{x \rightarrow \infty}
\frac{s_{B}(x)}{s_{A}(x - 2d)} = \alpha, \hspace*{3mm} \limsup_{x \rightarrow
  \infty} \frac{s_{B}(x)}{s_{A}(x + 2d)} = \beta.
\]
\end{enumerate}
\end{thm}
\noindent The next theorem shows that for a finite $s_{A}$ and finite or
infinite $s_{B}$ the rate of the unbounded distance $d_{A,B}(n)$ can be
arbitrary small.
\begin{thm}
Let $a \in \mathbb{Z^{+}}$, $b \in \mathbb{Z^{+}} \cup \{\infty\}$ and let $f:
\mathbb{Z^{+}} \rightarrow \mathbb{N}$ be a monoton increasing function
with $f(n) \rightarrow \infty$. Then there exist $A, B \subset \mathbb{N}$ such
that $s_{A} = a$, $s_{B} = b$ and $|a_{n} - b_{n}| \le f(n)$. The last
inequality implies, that $d_{A, B}(n-1) \le f(n)$.
\end{thm}
\noindent The right hand side of inequality (1) in Lemma 1. implies that for a
function $f(n)$ and for Sidon set $A$
and set $B$ having condition $d_{A,B}(n)\ll f(n)$ we have $ s_{B}(x) \ll
f(x)$. The next theorem tells us that this is sharp for $f(n)=n^{1/3}$.
\begin{thm}
There exist a Sidon set (i.e., $s_{A} = 1$) $A \subset \mathbb{N}$, $A =
\{a_{1}, a_{2}, \dots{}\}$ and a set $B \subset \mathbb{N}$, $B = \{b_{1},
b_{2}, \dots{} \}$ such that $d_{A,B}(n) \ll n^{1/3}$ and $s_{B}(n) \gg
 n^{1/3}$.
\end{thm}
\noindent If $s_A$ is finite then a routine calculation gives that $A(n)=O(\sqrt{n})$. Hence for $d_{A,B}(n)=O(n)$ we get that $B(n)=O(\sqrt{n})$, which implies that $s_B(n)=O(\sqrt{n})$. We pose the following question.

\begin{prob}
Is it true that for every $1/3 < \alpha < 1/2$ there exist sets $A, B
 \subset \mathbb{N}$ such that $s_{A} < \infty$, $d_{A,B}(n) \ll
 n^{\alpha}$ and $s_{B}(n) \gg n^{\alpha}$?
\end{prob}
\noindent The following theorem is about the case when both the representation
function $R_A(n)$ and the distance $d_{A,B}(n)$ are unbounded.
\begin{thm}
Let $A, B \subset \mathbb{N}$ such that $s_{A} = \infty$, $d_{A,B} = \infty$.
Then we have
\[
\max\{1, \frac{s_{A}(x - 2d_{A,B}(x))}{d_{A,B}(x)}\} \ll s_{B}(x) \ll
 s_{A}(x+2d_{A,B}(x))d_{A,B}(x).
\]
\end{thm}
\noindent We will prove Theorem 2. from Theorem 7. Starting out from the set of squares $A=\{n^2: n\geq 1 \}$ we will show that for any
$\varepsilon > 0$ arbitrary small, we have
$s_{A}(x)>exp\left( (\log 2-\varepsilon )\frac{\log x}{\log \log x}\right)$ for $x\geq x_0(\varepsilon)$. This implies Theorem 2.
The following three theorems deals with cubes.
\begin{thm}
Let $C$ be the set of positive cubes. Then we have $s_C(x)\gg \log \log x$.
\end{thm}
\noindent The next result is a direct consequence of Theorem 7 and Theorem 8.
\begin{cor}
Let us suppose that for a sequence $a_n$ we have $a_n=n^3+o(\log \log
 n)$. Then $s_A=\infty$
\end{cor}
\noindent In the other direction we have
\begin{thm}
There exists a set of positive integers such that $a_{n} = n^{3} +
O(n^{2.5}\log n)$, and $s_{A} < \infty$.
\end{thm}
\begin{prob} What conditions on $\alpha$ are needed to ensure that there
 exists a sequence $a_{n}$ such that $a_{n} = n^{3} +
O(n^{\alpha})$, and $s_{A} < \infty$?
\end{prob}
\noindent It follows from Theorem 7 that
\begin{cor}
Let $s_{A}(x) \sim x^{\alpha}$, $d_{A,B}(x) \sim x^{\beta}$. Then we have
\[
x^{\max\{0, \alpha - \beta\}} \ll s_{B}(x) \ll x^{\min\{\alpha + \beta, 1\}}.
\]
\end{cor}

\begin{prob}
Let $\alpha$, $\beta$ be nonnegative real numbers and max $\{0, \alpha -
\beta\} \le \gamma \le \min \{\alpha + \beta, 1\}$. What conditions on
$\alpha$, $\beta$ and $\gamma $ are needed to ensure that there exist
$A$, $B \subset \mathbb{N}$ such that $s_{A}(x) \sim x^{\alpha}$, $d_{A,B}(x) \sim x^{\beta}$ and $s_{B}(x) \sim x^{\gamma}$?
\end{prob}

\section{Proofs}

\begin{lem} For any subsets $A$ and $B$ of $\mathbb{N}$ we have
\begin{equation}
\frac{s_{A}(x-2d_{A,B}(x))}{4d_{A,B}(x)+1} \le s_{B}(x) \le
 s_{A}(x+2d_{A,B}(x))(4d_{A,B}(x)+1).
\end{equation}
\end{lem}
\noindent {\bf Proof.} Let $n \le x$. If $b_{i} + b_{j} = n$, then $b_{i} \le
n \le x$ and $b_{j} \le n \le x$. By the definition of $d_{A,B}(x)$ we have
$$
|a_{i} - b_{i}| \le \max_{\hbox{t: $a_{t} \le b_{i}$ or $b_{t} \le b_{i}$}}|a_{t} - b_{t}| = d_{A,B}(b_{i}),
$$
then we have
\[
-d_{A,B}(b_{i}) \le a_{i} - b_{i} \le d_{A,B}(b_{i}),
\]
thus
\[
b_{i} - d_{A,B}(b_{i}) \le a_{i} \le d_{A,B}(b_{i}) + b_{i},
\]
which implies
\[
b_{i} - d_{A,B}(x) \le a_{i} \le d_{A,B}(x) + b_{i}.
\]
Similarly for $b_{j}$,
\[
b_{j} - d_{A,B}(x) \le a_{j} \le d_{A,B}(x) + b_{j}.
\]
It follows that
\[
b_{i} + b_{j} - 2d_{A,B}(x) \le a_{i} + a_{j} \le 2d_{A,B}(x) + b_{j} + b_{i},
\]
thus
\[
n - 2d_{A,B}(x) \le a_{i} + a_{j} \le 2d_{A,B}(x) + n.
\]
Then we have
\[
R_{B}(n) \le |\{(i,j): n - 2d_{A,B}(x) \le a_{i} + a_{j} \le n +
2d_{A,B}(x), i \le j\}|
\]
\[ \displaystyle
\le \sum_{n - 2d_{A,B}(x) \le m \le n+2d_{A,B}(x)}R_{A}(m) \le (4d_{A,B}(x)+1)
\max_{n - 2d_{A,B}(x) \le m \le n+2d_{A,B}(x)}R_{A}(m)
\]
\[
\le s_{A}(x+2d_{A,B}(x))(4d_{A,B}(x)+1),
\]
which proves the second inequality of (1). If we replace $x$ by $x -
2d_{A,B}(x)$, and $A$ by $B$ in the second inequality and by using the
fact that $d_{A,B}(x-2d_{A,B}(x)) \le d_{A,B}(x)$, we obtain that
\[
s_{A}(x-2d_{A,B}(x)) \le  s_{B}(x-2d_{A,B}(x) +
2d_{A,B}(x-2d_{A,B}(x)))(4d_{A,B}(x-2d_{A,B}(x))+1))
\]
\[
\le s_{B}(x)(4d_{A,B}(x)+1),
\]
which proves the first inequality in (1). $\blacksquare $

\bigskip

\noindent {\bf Proof of the first part of Theorem 3.}
It is clear that there exists an $x_{0}$ such that if $x \ge x_{0}$ then
$s_{A}(x) = s_{A}$, $s_{B}(x) = s_{B}$, $d_{A,B}(x) = d_{A,B}$. If we
choose $x = x_{0} + 2d_{A,B}$ it follows from (1) that
\[
\frac{s_{A}(x_{0}+2d_{A,B}-2d_{A,B}(x))}{4d_{A,B}(x_{0} + 2d_{A,B})+1}
\le s_{B}(x_{0} + 2d_{A,B}) \le s_{A}(x_{0}+4d_{A,B}(x))(4d_{A,B}(x_{0}
+ 2d_{A,B})+1),
\]
thus
\[
\frac{s_{A}}{4d_{A,B}+1} \le s_{B} \le (4d_{A,B}+1)s_{A},
\]
which proves the first part of Theorem 3. $\blacksquare$

\bigskip

\noindent {\bf Proof of the first part of Theorem 4.}
The first and the third inequality follows from Lemma 1 immediately. We
prove the second inequality by contradiction. Assume that
\[
\liminf_{x\rightarrow \infty}\frac{s_{B}(x)}{s_{A}(x-2d_{A,B})} =
\alpha, \hspace*{5mm} \limsup_{x\rightarrow
\infty}\frac{s_{B}(x)}{s_{A}(x+2d_{A,B})} = \beta,
\]
where $\alpha > \beta$. Then for any $\varepsilon > 0$ there exist an
$x_{0}$ such that for $x \ge x_{0}$ we have $s_{B}(x) > (\alpha -
\varepsilon)s_{A}(x-2d_{A,B})$ and $s_{B}(x) < (\beta +
\varepsilon)s_{A}(x+2d_{A,B})$. Let $x_{0} \le N \le M$. Then we have
\[
\sum_{N \le x \le M}(\alpha - \varepsilon)s_{A}(x-2d_{A,B}) < \sum_{N
\le x \le M}(\beta + \varepsilon)s_{A}(x+2d_{A,B}),
\]
thus
\begin{equation}
\sum_{N+2d_{A,B} \le x \le M-2d_{A,B}}(\alpha - \beta - 2\varepsilon)s_{A}(x)
\end{equation}
\[
\le \sum_{M-2d_{A,B} < x \le M+2d_{A,B}}(\beta + \varepsilon)s_{A}(x)
\]
\[
\le 4(\beta + \varepsilon)s_{A}(M+2d_{A,B})d_{A,B}.
\]
Let $\varepsilon = \frac{\alpha - \beta}{3}$. We may choose $N$ such that
\[
s_{A}(N) \ge \frac{4(\beta + \varepsilon)4d_{A,B}}{\alpha -
\beta}.
\]
Then we have
\[
\sum_{N+2d_{A,B} < x \le M-2d_{A,B}}(\alpha - \beta -
2\varepsilon)s_{A}(x) \ge \frac{\alpha - \beta}{3}(M - N -
4d_{A,B})s_{A}(N) \ge
\]
\[
\frac{\alpha - \beta}{3} \cdot \frac{4(\beta +
\varepsilon)4d_{A,B}}{\alpha - \beta}(M - N - 4d_{A,B}) =
\frac{4}{3}(\beta + \varepsilon)4(M - N - 4d_{A,B})d_{A,B}.
\]
It follows from (2) that
\[
\frac{4}{3}(\beta + \varepsilon)4(M - N - 4d_{A,B})d_{A,B} \le 4(\beta +
\varepsilon)(M+2d_{A,B})d_{A,B},
\]
thus we have
\[
\frac{4}{3}(M - N - 4d_{A,B}) \le M+2d_{A,B},
\]
which is a contradiction if $N$ is fixed and $M$ is large
enough. $\blacksquare$

\bigskip

\noindent {\bf Proof of Theorem 7.}
Theorem 7. follows from Lemma 1 immediately. $\blacksquare$

\bigskip

\noindent {\bf Proof of Theorem 2.}
Let $A$ be the set of positive squares. Consider the sequence primes $q_1 < q_2 < \dots{}$ where every $q_i \equiv 1 \bmod{4}$. Define
\[
Q_{k} = \prod_{i=1}^{k}q_{i}.
\]
Let $Q_{K}$ denote the maximal $Q_{k}$ not exceeding $x$. It is easy to see from the well known formula about the number of representations of a positive integer as the sum of two squares [23] that $R_{A}(Q_{K}) = 2^{K-1}$. The well known formula for primorial [27] implies that
\[
Q_{K} = e^{(1+o(1))K\log K}.
\]
It follows that
\[
R_{A}(Q_{K}) = e^{(\log 2+o(1))\frac{\log Q_{K}}{\log \log Q_{K}}}.
\]
It is clear that
\[
\frac{\log (Q_{k+1})}{\log (Q_{k})} \rightarrow 1
\]
as $k \rightarrow \infty$. Thus we have
\[
s_{A}(x) \ge e^{(\log 2+o(1))\frac{\log x}{\log \log x}}.
\]
We apply Theorem 7. If $b_{n} = n^{2} + O\Big(e^{(\log 2-\varepsilon)\frac{\log x}{\log \log x}}\Big)$, then $d_{A,B}(x) = O\Big(e^{(\log 2-\varepsilon)\frac{\log x}{\log \log x}}\Big)$. Then by Theorem 7. we have
\[
s_{B}(x) \gg \frac{s_{A}(x - 2d_{A,B}(x))}{d_{A,B}(x)} \ge \frac{e^{(\log 2+o(1))\frac{\log x}{\log \log x}}}{e^{(\log 2-\varepsilon)\frac{\log x}{\log \log x}}} = e^{(\varepsilon+o(1))\frac{\log x}{\log \log x}},
\]
which implies that $s_{B} = \infty$, thus the function $R_{B}(n)$ is unbounded.
$\blacksquare$

\bigskip

\noindent To prove the second part of Theorem 3. and Theorem 4. we need
the following Lemma.

\begin{lem} Let $a(n)$, $b(n)$, $d(n)$ be monotone increasing sequences of
 positive integers satisfying
\[
\frac{a(n)}{4d(n)+1} \le b(n) \le a(n)(4d(n)+1).
\]
We put
\[
T(n) = \sum_{k \le n}max\{a(k), b(k)\}.
\]
Let $c_{n} = d(n)(2 + 10^{n+T(n)})$. Then there exist sets $A, B \subset
\mathbb{N}$ such that the following holds.
\begin{itemize}
\item[(i)] $s_{A}(c_{n} - 2d(n)) = s_{A}(c_{n}+2d(n)) = a(n)$,
\item[(ii)] $s_{B}(c_{n} - 2d(n)) = s_{B}(c_{n}+2d(n)) = b(n)$,
\item[(iii)] $d_{A,B}(c_{n}) = d(n)$,
\item[(iv)] $A(c_{n}) = B(c_{n}) = 2T(n) \ge 2n$.
\end{itemize}
\end{lem}

\noindent {\bf Proof.} We give a construction for the sets $A$, $B$,
which built up by blocks.
Let $\displaystyle A = \cup_{n=1}^{\infty}A^{(n)}$, and $\displaystyle B
= \cup_{n=1}^{\infty}B^{(n)}$, where
\[
A^{(n)} = \{a^{(n)}_{i}: 1 \le i \le 2max\{a(n), b(n)\}\}.
\]
\[
B^{(n)} = \{b^{(n)}_{i}: 1 \le i \le 2max\{a(n), b(n)\}\}.
\]
Assume that $max\{a(n), b(n)\} = a(n)$. Then let
\[
a^{(n)}_{i} =
\begin{cases}
d(n)10^{n - 1 + i + T(n-1)}, & \text{if $1 \le i \le a(n)$}\\
d(n)10^{n + T(n)} - d(n)10^{n - 1 + i - a(n) + T(n-1)}, & \text{if
 $a(n)+1 \le i \le 2a(n)$},
\end{cases}
\]
therefore $a^{(n)}_{i}=c_{n} -
 2d(n) - a^{(n)}_{i-a(n)}$ for $a(n)+1 \le i \le 2a(n)$. It is easy to see that for any $1 \le i \le a(m)$ and $a \in A$, $a < a^{(m)}_{i}$ we have
\begin{equation}
10a \le a^{(m)}_{i}.
\end{equation}
Let
\[
b^{(n)}_{i} =
\begin{cases}
d(n)10^{n - 1 + i + T(n-1)} - d(n) - 1 + \Big\lceil\frac{i}{2b(n)}\Big\rceil,
& \text{if $1 \le i \le a(n)$}\\
d(n)10^{n + T(n)} - d(n)10^{n - 1 + i - a(n) + T(n)} \\ \qquad {}- d(n) -
 1 + \Big\lceil \frac{i - a(n) +
 b(n)}{2b(n)}\Big\rceil,
& \text{if $a(n)+1 \le i \le 2a(n)$}.
\end{cases}
\]
Assume that $max\{a(n), b(n)\} = b(n)$. Then let
\[
a^{(n)}_{i} =
\begin{cases}
d(n)10^{n - 1 + i + T(n-1)} - d(n) - 1 + \Big\lceil
 \frac{i}{2a(n)}\Big\rceil, & \text{if $1 \le i \le b(n)$}\\
d(n)10^{n + T(n)} - d(n)10^{n - 1 + i - b(n) + T(n-1)} - d(n) - 1 + \Big\lceil
 \frac{i - b(n) + a(n)}{2a(n)}\Big\rceil, & \text{if $b(n)+1 \le i
 \le 2b(n)$},
\end{cases}
\]
and
\[
b^{(n)}_{i} =
\begin{cases}
d(n)10^{n - 1 + i + T(n-1)},
& \text{if $1 \le i \le b(n)$}\\
d(n)10^{n + T(n)} - d(n)10^{n - 1 + i - b(n) + T(n-1)}, & \text{if
 $b(n)+1 \le i \le 2b(n)$}.
\end{cases}
\]
 It is easy to see that for any $1 \le i \le b(m)$ and $a \in A$, $a < a^{(m)}_{i}$ we have
\begin{equation}
10a \le a^{(m)}_{i}.
\end{equation}

In the next step we prove that the above sets $A$ and $B$ satisfy the
Lemma.

First we prove $(i)$. To do this we show that both $s_{A}(c_{n} + 2d(n)) \ge a(n)$ and $s_{A}(c_{n} + 2d(n)) \le a(n)$ hold. In the first case assume that $max\{a(n), b(n)\} = a(n)$. It is clear from the definition of $a^{(n)}_{i}$ that $a^{(n)}_{i} + a^{(n)}_{i+a(n)}  = d(n)10^{n + T(n)} = c_{n} - 2d(n)$ for $1 \le i \le a(n)$. Therefore, $R_{A}\Big(d(n)10^{n + T(n)}\Big) = R_{A}(c_{n} - 2d(n)) \geq a(n)$. Thus we have $s_{A}(c_{n} + 2d(n)) \ge a(n)$. In the second case we assume that $max\{a(n), b(n)\} = b(n)$. Now we have $a^{(n)}_{i} + a^{(n)}_{i+b(n)}  = d(n)10^{n + T(n)} - 2d(n) = c_{n} - 4d(n)$ for $1 \le i \le a(n)$, which implies that $R_{A}(c_{n} - 4d(n)) \ge a(n)$ and so $s_{A}(c_{n} + 2d(n)) \ge a(n)$. It remains to prove that  $s_{A}(c_{n} + 2d(n)) \le a(n)$. Since $d(n)$ is monotone increasing, if $max\{a(n), b(n)\} = a(n)$ then $a^{(n+1)}_{1} = d(n+1)10^{n + 1 + T(n)} > c_{n} + 2d(n) = d(n)10^{n + T(n)} + 4d(n)$. On the other hand, if $max\{a(n), b(n)\} = b(n)$ then $a^{(n+1)}_{1} = d(n+1)10^{n + 1 + T(n+1)} - d(n+1) > c_{n} + 2d(n) = d(n)10^{n + T(n)} + 4d(n)$. It is enough to show that $R_{A}(a_{s} + a_{t}) \le a(n)$ for any fixed
\[
a_{s}, a_{t} \in \cup_{i=1}^{n}A^{(i)}.
\]
In the next we give an upper estimation to the number of pairs
$(u, v)$ such that $a_s+a_t=a_u+a_v$, where
$a_u,a_v\in A$, and we may suppose that $a_{s}\geq a_u\geq a_v\geq a_t$, and $a_{s} \in A^{(m)}$ for some $m \le n$.
We distinguish two cases.

In the first case assume that $max\{a(m), b(m)\} = a(m)$. We will prove that if $a_{s} + a_{t} = a_{u}
+ a_{v}$ is a nontrivial solution, i.e., $R_{A}(a_{s} + a_{t}) \ge 2$, then
\[
a_{s} + a_{t} = d(m)10^{m + T(m)}.
\]
We have five subcases depending on how many of $a_{s}, a_{t}$, $a_{u}$,
$a_{v}$ are selected from the set $\{a_{1+10^m}^{(m)},a_{2+10^m}^{(m)},\dots ,a_{10^m+10^m}^{(m)} \}$.

If none of them are selected from the set $\{a_{1+10^m}^{(m)},a_{2+10^m}^{(m)},\dots ,a_{10^m+10^m}^{(m)} \}$, then $a_{j}^{(m)} + a_{t} =
a_{u} + a_{v}$. It follows from (3) that $a_{u} = a_{j}^{(m)}$, thus
$a_{s} + a_{t} = a_{u} + a_{v}$ must be a trivial solution.

If one of the terms $a_{s}$, $a_{t}$, $a_{u}$, $a_{v}$ are selected from the set $\{a_{1+10^m}^{(m)},a_{2+10^m}^{(m)},\dots ,a_{10^m+10^m}^{(m)} \}$, that is, $a_{s} = a_{j}^{(m)}$, for some $a(m) < j \le 2a(m)$ then it is clear that $a_{j}^{(m)} + a_{t} > a_{u} + a_{v}$, a contradiction.

If two of the terms $a_{s}$, $a_{t}$, $a_{u}$, $a_{v}$ are selected from the set $\{a_{1+10^m}^{(m)},a_{2+10^m}^{(m)},\dots ,a_{10^m+10^m}^{(m)} \}$, then let $a_{s} = a_{j}^{(m)}$ and $a_{u} = a_{k}^{(m)}$, where $a(m) < j \le k \le 2a(m)$. Then we have
\[
c_{m} - 2d(m) - a_{j-a(m)}^{(m)} + a_{t} = c_{m} - 2d(m) - a_{k-a(m)}^{(m)}
+ a_{v}.
\]
Thus we have

\begin{equation}
a_{k-a(m)}^{(m)} + a_{t} = a_{j-a(m)}^{(m)} + a_{v}.
\end{equation}
If $k=j$ then $a_s+a_t=a_u+a_v$ is a trivial solution. If $j<k$ then it follows from (3) that $a_v=a_{k-a(m)}^{(m)}$ and
$a_t=a_{j-a(m)}^{(m)}$, therefore
$a_s+a_t=a_j^{(m)}+a_{j-a(m)}^{(m)}=d(m)10^{m+T(m)}$.

If three of the terms $a_{s}$, $a_{t}$, $a_{u}$, $a_{v}$ are selected
from the set $\{a_{1+10^m}^{(m)},a_{2+10^m}^{(m)},\dots
,a_{10^m+10^m}^{(m)} \}$, then let $a_{s} = a_{j}^{(m)}$ and
$a_u=a_{k}^{(m)}$, where $a(m) < j \le k \le 2a(m)$. Then we have
$$a_{k-a(m)}^{(m)}+a_t=a_{j-a(m)}^{(m)}+a_v,$$
where only one term is selected from the set
$\{a_{1+10^m}^{(m)},a_{2+10^m}^{(m)},\dots ,a_{10^m+10^m}^{(m)} \}$,
which is absurd.

If four of the terms $a_{s}$, $a_{t}$, $a_{u}$, $a_{v}$ are selected
from the set $\{a_{1+10^m}^{(m)},a_{2+10^m}^{(m)},\dots
,a_{10^m+10^m}^{(m)} \}$, then let $a_{s} = a_{j}^{(m)}$,
$a_u=a_{k}^{(m)}$, $a_v=a_{l}^{(m)}$ and $a_t=a_{q}^{(m)}$ where $a(m) <
j \le k \le l \le q \le 2a(m)$. Then
$$a_{k-a(m)}^{(m)}+a_{l-a(m)}^{(m)}=a_{j-a(m)}^{(m)}+a_{q-a(m)}^{(m)},$$
 which must be a trivial solution.

We have $a(m)$ elements of the set $A$ in the interval
$$\Big[\frac{d(m)}{2}10^{m+T(m)}, d(m)10^{m+T(m)}\Big],$$ therefore $R_{A}(a_{s} + a_{t}) \le a(m) \le a(n)$.

In the second case assume that $max\{a(m), b(m)\} = b(m)$. In this case we prove that the nontrivial equation $a_{s} + a_{t} = a_{u} + a_{v}$ implies the existence of the integers $1 \le j \le k \le b(m)$ such that $a_{s} = a_{j}^{(m)}$, $a_{t} =  a_{j-b(m)}^{(m)}$, $a_{u} = a_{k}^{(m)} $ and $a_{v} =  a_{k-b(m)}^{(m)}$.

Then we have five subcases as above.

If none of them are selected from the set $\{a_{1+10^m}^{(m)},a_{2+10^m}^{(m)},\dots ,a_{10^m+10^m}^{(m)} \}$, then $a_{s} + a_{t} =
a_{u} + a_{v}$, where $a_{s} = a_{j}^{(m)}$ and $1 \le j \le b(m)$. It follows from (4) that $a_{u} = a_{j}^{(m)}$, thus
$a_{s} + a_{t} = a_{u} + a_{v}$ must be a trivial solution.

If one of the terms $a_{s}$, $a_{t}$, $a_{u}$, $a_{v}$ are selected from the set $\{a_{1+10^m}^{(m)},a_{2+10^m}^{(m)},\dots ,a_{10^m+10^m}^{(m)} \}$, that is, $a_{s} = a_{j}^{(m)}$, for some $b(m) < j \le 2b(m)$ then it is clear that $a_{j}^{(m)} + a_{t} > a_{u} + a_{v}$, a contradiction.

If two of the terms $a_{s}$, $a_{t}$, $a_{u}$, $a_{v}$ are selected from the set $\{a_{1+10^m}^{(m)},a_{2+10^m}^{(m)},\dots ,a_{10^m+10^m}^{(m)} \}$, then let $a_{s} = a_{j}^{(m)}$ and $a_{u} = a_{k}^{(m)}$, where $b(m) < j \le k \le 2b(m)$. Then the equation $a_{s} + a_{t} = a_{u} + a_{v}$ means that
$$d(m)10^{m + T(m)} - d(m)10^{m - 1 + j - b(m) + T(m-1)} - d(m)
- 1 + \Big\lceil \frac{j-b(m)+a(m)}{2a(m)}\Big\rceil + a_{t}$$
$$ = d(m)10^{m + T(m)} - d(m)10^{m - 1 + k - b(m) + T(m-1)} - d(m)
- 1 + \Big\lceil \frac{k-b(m)+a(m)}{2a(m)}\Big\rceil + a_{v},$$
that is
$$ d(m)10^{m - 1 + k - b(m) + T(m-1)} + a_{t}$$
$$ = d(m)10^{m - 1 + j - b(m) + T(m-1)} + a_{v} + \Big\lceil \frac{k-b(m)+a(m)}{2a(m)}\Big\rceil - \Big\lceil \frac{j-b(m)+a(m)}{2a(m)}\Big\rceil,$$
therefore,
$$d(m)10^{m - 1 + k - b(m) + T(m-1)} - d(m) - 1 + \Big\lceil \frac{k-b(m)}{2a(m)}\Big\rceil + a_{t}$$
$$ = d(m)10^{m - 1 + j - b(m) + T(m-1)} - d(m) - 1 + \Big\lceil \frac{j-b(m)}{2a(m)}\Big\rceil + a_{v}$$
$$ + \Big\lceil \frac{k-b(m)}{2a(m)}\Big\rceil - \Big\lceil \frac{j-b(m)}{2a(m)}\Big\rceil + \Big\lceil \frac{k-b(m)+a(m)}{2a(m)}\Big\rceil -
\Big\lceil \frac{j-b(m)+a(m)}{2a(m)}\Big\rceil,$$
i.e.,
$$a_{k-b(m)}^{(m)}+a_{t}=a_{j-b(m)}^{(m)} + a_{v} + \Big\lceil \frac{k-b(m)}{2a(m)}\Big\rceil + \Big\lceil \frac{k-b(m)+a(m)}{2a(m)}\Big\rceil $$
$$-\Big(\Big\lceil \frac{j-b(m)}{2a(m)}\Big\rceil + \Big\lceil \frac{j-b(m)+a(m)}{2a(m)}\Big\rceil\Big).$$
If $k=j$ then $a_s+a_t=a_u+a_v$ is a trivial solution. Since $b(m) \le a(m)(4d(m)+ 1)$ and $b(m) < j \le k \le 2b(m)$, if $j<k$ then we have
\[
1 \le \Big\lceil \frac{k - b(m)}{2a(m)}\Big\rceil \le \Big\lceil
\frac{b(m)}{2a(m)}\Big\rceil \le \Big\lceil \frac{(4d(m)+1)a(m)}{2a(m)}\Big\rceil = 2d(m) + 1,
\]
and
\[
1 \le \Big\lceil \frac{k - b(m) + a(m)}{2a(m)}\Big\rceil \le \Big\lceil \frac{(4d(m)+1)a(m) + a(m)}{2a(m)}\Big\rceil = 2d(m) + 1.
\]
It follows from (4) that $a_v=a_{k-b(m)}^{(m)}$ and by the definition of $a_{j-b(m)}^{(m)} = d(m)10^{m - 1 + j - b(m) + T(m-1)} $ it follows that $a_t=a_{j-b(m)}^{(m)}$. Fix the elements $a_{s}$ and $a_{t}$. The equation $a_s+a_t=a_u+a_v$ can be written in the form $a_{j}^{(m)} + a_{j-b(m)}^{(m)} = a_{k}^{(m)} + a_{k-b(m)}^{(m)}$.

If three of the terms $a_{s}$, $a_{t}$, $a_{u}$, $a_{v}$ are selected
from the set $\{a_{1+10^m}^{(m)},a_{2+10^m}^{(m)},\dots
,a_{10^m+10^m}^{(m)} \}$, then let $a_{s} = a_{j}^{(m)}$,
$a_u=a_{k}^{(m)}$ and $a_v=a_{l}^{(m)}$ , where $b(m) < j \le k \le l \le 2b(m)$. Then we have
$$a_{j}^{(m)}+a_t < a_{k}^{(m)}+a_{l}^{(m)},$$
where only one term is selected from the set
$\{a_{1+10^m}^{(m)},a_{2+10^m}^{(m)},\dots ,a_{10^m+10^m}^{(m)} \}$,
which is absurd as we have seen earlier.

If four of the terms $a_{s}$, $a_{t}$, $a_{u}$, $a_{v}$ are selected
from the set $\{a_{1+10^m}^{(m)},a_{2+10^m}^{(m)},\dots
,a_{10^m+10^m}^{(m)} \}$, then let $a_{s} = a_{j}^{(m)}$,
$a_u=a_{k}^{(m)}$, $a_v=a_{l}^{(m)}$ and $a_t=a_{m}^{(m)}$ where $b(m) <
j \le k \le l \le q \le 2b(m)$. Then by definition the equation
$a_{j}^{(m)}+a_{q}^{(m)}=a_{k}^{(m)}+a_{l}^{(m)}$ means that
$$d(m)10^{m + T(m)} - d(m)10^{m - 1 + j - b(m) + T(m-1)} - d(m)
- 1 + \Big\lceil \frac{j-b(m)+a(m)}{2a(m)}\Big\rceil $$
$$ + d(m)10^{m + T(m)} - d(m)10^{m - 1 + q - b(m) + T(m-1)} - d(m)
- 1 + \Big\lceil \frac{q-b(m)+a(m)}{2a(m)}\Big\rceil $$
$$ = d(m)10^{m + T(m)} - d(m)10^{m - 1 + k - b(m) + T(m-1)} - d(m)
- 1 + \Big\lceil \frac{k-b(m)+a(m)}{2a(m)}\Big\rceil $$
$$ + d(m)10^{m + T(m)} - d(m)10^{m - 1 + l - b(m) + T(m-1)} - d(m)
- 1 + \Big\lceil \frac{l-b(m)+a(m)}{2a(m)}\Big\rceil. $$
Therefore,
$$ d(m)10^{m - 1 + j - b(m) + T(m-1)} - \Big\lceil \frac{j-b(m)+a(m)}{2a(m)}\Big\rceil + d(m)10^{m - 1 + q - b(m) + T(m-1)} - \Big\lceil \frac{q-b(m)+a(m)}{2a(m)}\Big\rceil$$
$$ = d(m)10^{m - 1 + k - b(m) + T(m-1)} - \Big\lceil \frac{k-b(m)+a(m)}{2a(m)}\Big\rceil + d(m)10^{m - 1 + l - b(m) + T(m-1)} - \Big\lceil \frac{l-b(m)+a(m)}{2a(m)}\Big\rceil, $$
which implies that $q = l$ and $j = k$, a trivial solution.

If $R_{A}(a_{s}+a_{t}) > 1$ then for the nontrivial solution $a_{s} + a_{t} = a_{u} + a_{v}$ we have integers $1 \le j \le k \le b(m)$ such that $a_{s} = a_{j}^{(m)}$ $a_{t} = a_{j-b(m)}^{(m)}$, $a_{u} = a_{k}^{(m)}$, $a_{v} = a_{k-b(m)}^{(m)}$.

By the definition we have
$$d(m)10^{m + T(m)} - d(m)10^{m - 1 + j - b(m) + T(m-1)} - d(m) - 1 + \Big\lceil
 \frac{j - b(m) + a(m)}{2a(m)}\Big\rceil $$
$$ + d(m)10^{m - 1 + j - b(m) + T(m-1)} - d(m) - 1 + \Big\lceil
 \frac{j - b(m)}{2a(m)}\Big\rceil $$
$$ = d(m)10^{m + T(m)} - d(m)10^{m - 1 + k - b(m) + T(m-1)} - d(m) - 1 + \Big\lceil
 \frac{k - b(m) + a(m)}{2a(m)}\Big\rceil $$
$$ + d(m)10^{m - 1 + k - b(m) + T(m-1)} - d(m) - 1 + \Big\lceil
 \frac{k - b(m)}{2a(m)}\Big\rceil. $$
Thus we have
$$\Big\lceil \frac{k - b(m) + a(m)}{2a(m)}\Big\rceil + \Big\lceil \frac{k - b(m)}{2a(m)}\Big\rceil = \Big\lceil \frac{j - b(m) + a(m)}{2a(m)}\Big\rceil + \Big\lceil \frac{j - b(m)}{2a(m)}\Big\rceil. $$
Let $j = (l_{1} - 2)a(m) + h_{1}$ and $k = (l_{2} - 2)a(m) + h_{2}$, where $1 \le h_{1}, h_{2} \le a(m)$.
Then we have
$$ \Big\lceil \frac{(l_{1}-2)a(m)+h_{1}}{2a(m)}\Big\rceil + \Big\lceil \frac{(l_{1}-1)a(m)+h_{1}}{2a(m)}\Big\rceil =
\Big\lceil \frac{(l_{2}-2)a(m)+h_{2}}{2a(m)}\Big\rceil +
\Big\lceil \frac{(l_{2}-1)a(m)+h_{2}}{2a(m)}\Big\rceil,
$$ hence
$$\Big\lceil \frac{l_{1}-1}{2}\Big\rceil +
\Big\lceil \frac{l_{1}}{2}\Big\rceil = \Big\lceil \frac{l_{2}-1}{2}\Big\rceil +\Big\lceil \frac{l_{2}}{2}\Big\rceil,$$
that is $l_{1} = l_{2}$, which implies that the number of suitable pairs $(u, v)$ is at most $a(m) \le a(n)$.

Changing the role of $a(n)$ and $b(n)$, the proof of (ii) is the same as the proof (i) thus we omit it and leave the details to the reader.

Now we prove (iii). Assume that $max\{a(m), b(m)\} = a(m)$, where $m \le n$. Then if
$1 \le i \le a(m)$, then
\[
a_{i}^{(m)} - b_{i}^{(m)} = d(m) + 1 - \Big\lceil \frac{i}{2b(m)}\Big\rceil.
\]
As $1 \le i \le a(m) \le (4d(m)+1)b(m)$, then we have
\[
1 \le \Big\lceil \frac{i}{2b(m)}\Big\rceil \le \Big\lceil
\frac{(4d(m)+1)b(m)}{2b(m)}\Big\rceil = 2d(m) + 1,
\]
then $-d(m) \le a_{i}^{(m)} - b_{i}^{(m)} \le d(m)$. If $a(m) < i \le
2a(m)$, then
\[
a_{i}^{(m)} - b_{i}^{(m)} = d(m) + 1 - \Big\lceil
\frac{i-a(m)+b(m)}{2b(m)}\Big\rceil.
\]
Since $a(m) \le (4d(m)+1)b(m)$, then
\[
1 \le \Big\lceil \frac{i - a(m) + b(m)}{2b(m)}\Big\rceil \le \Big\lceil
\frac{a(m) + b(m)}{2b(m)}\Big\rceil \le \Big\lceil \frac{(4d(m)+1)b(m) +
b(m)}{2b(m)}\Big\rceil = 2d(m) + 1,
\]
thus $-d(m) \le a_{i}^{(m)} - b_{i}^{(m)} \le d(m)$. As $d(m)$ is
monotonous increasing and $|a_{1}^{(m)}-b_{1}^{(m)}|=d(m)$, it follows (iii).
If $max\{a(m), b(m)\} = b(m)$, changing the role of $a(m)$ and $b(m)$, the argument is the same.

The statement (iv) follows easily from the construction. By the definition of $A^{(n)}$    it is clear that it has at least two elements. Therefore, up to $c_{n}$ the set $A$ has at least $2n$ elements. The same argument works for the set $B$ as well. $\blacksquare $

\bigskip

\noindent {\bf Proof of the second part of Theorem 3.}
If we put $a(n) = a$, $b(n) = b$, $d(n) = d$ in Lemma 2, the Theorem follows
immediately. $\blacksquare$

\bigskip

\noindent {\bf Proof of the second part of Theorem 4.}
Let us suppose that the monotone increasing sequences $\{u_{n}\}$ and
$\{v_{n}\}$ satisfy
\[
\lim_{n \rightarrow \infty}\frac{u_{n}}{v_{n}} = \alpha,
\]
and
\[
\lim_{n \rightarrow \infty}\frac{u_{n+1}}{v_{n}} = \beta.
\]
Let $b(2n-1) = u_{n}$, $b(2n) = u_{n+1}$, $a(2n-1) = a(2n) = v_{n}$,
and $d(n) = d_{A,B}$ in Lemma 2. then by the construction we have
\[
\frac{s_{B}(c_{2n-1})}{s_{A}(c_{2n-1}-2d_{A,B})} = \frac{u_{n}}{v_{n}}
\rightarrow \alpha,
\]
and
\[
\frac{s_{B}(c_{2n})}{s_{A}(c_{2n}+2d_{A,B})} = \frac{u_{n+1}}{v_{n}}
\rightarrow \beta.
\]
If $c_{2n-1} \le k \le c_{2n}$, then $s_{A}(c_{2n-1}-2d_{A,B}) =
s_{A}(k-2d_{A,B}) = s_{A}(k+2d_{A,B}) =  s_{A}(c_{2n}+2d_{A,B}) = v_{n}$.
It follows that
\[
\alpha \leftarrow \frac{s_{B}(c_{2n-1})}{s_{A}(c_{2n-1}-2d_{A,B})} \le
\frac{s_{B}(k)}{s_{A}(k-2d_{A,B})} =
\frac{s_{B}(k)}{s_{A}(k+2d_{A,B})} \le
\frac{s_{B}(c_{2n})}{s_{A}(c_{2n} + 2d_{A,B})} \rightarrow \beta.
\]
If $c_{2n} \le k \le c_{2n+1}$, then $s_{B}(c_{2n}-2d_{A,B}) = s_B(k)
= s_{B}(c_{2n+1}+2d_{A,B}) = u_{n+1}$. It follows that
\[
\alpha \leftarrow \frac{s_{B}(c_{2n+1})}{s_{A}(c_{2n+1}-2d_{A,B})} \le
\frac{s_{B}(k)}{s_{A}(k+2d_{A,B})} \le
\frac{s_{B}(k)}{s_{A}(k-2d_{A,B})} \le
\frac{s_{B}(c_{2n})}{s_{A}(c_{2n}+2d_{A,B})} \rightarrow \beta,
\]
which completes the proof. $\blacksquare$

\bigskip

\noindent {\bf Proof of Theorem 5.}
First assume that $b$ is a positive integer. Let $a(n) = a$.
Without loss of generality we may assume that $a \le b$. Let
\[
b(n) =
\begin{cases}
a, & \text{if $f(n) < \Big\lceil\frac{b}{a}\Big\rceil$}\\ \\
b, & \text{if $f(n) \ge \Big\lceil\frac{b}{a}\Big\rceil$},
\end{cases}
\]
and let $d(n) = f(n)$ in Lemma 2. Since
$A(c_{n}) \ge 2n$, $d_{A,B}(c_{n}) = d(n) = f(n)$, then if $k \le
2n$, this implies that $|a_{k} - b_{k}| \le f(n)$, therefore $|a_{n} - b_{n}|
\le f(n)$.

On the other hand, if $b = +\infty$, then let $a(n) = a$, $b(n) = af(n)$,
$d(n) = f(n)$ in Lemma 2. It follows that $s_{B} = +\infty$, $s_{A} = a$,
$A(c_{n}) \ge 2n$, $d_{A,B}(c_{n}) = d(n) = f(n)$. This implies that
$|a_{n} - b_{n}| \le f(n)$. $\blacksquare$

\bigskip

\noindent {\bf Proof of Theorem 6.}
We give a construction for the sets $A$ and $B$ recursively. Define the
sets $A^{(m)}$ and
$B^{(m)}$ in the following way. Fix a nonnegative integer $m$ and we will
choose the distinct positive integers
$200 \cdot 1000^{m} \le b_{1}^{(m)}, \dots{} ,b_{10^{m}}^{(m)} \le 300 \cdot
1000^{m}$. Define $b_{i+10^m}^{(m)}$ by
$b_{i+10^m}^{(m)} =
1000^{m+1} - b_{i}^{(m)}$, where $1 \le i \le 10^{m}$. Let $B^{(m)} = \{b_{1}^{(m)}, \dots{} ,b_{10^{m}}^{(m)}, b_{10^m+1}^{(m)}, \dots{}
,b_{10^m+10^m}^{(m)}\}$. It is easy to see that $R_{B}(1000^{m+1}) \ge
10^{m} \gg (1000^{m+1})^{1/3}$. This implies that $s_{B}(n) \gg n^{1/3}$. Now we define the sets $A^{(m)}$ in the
following way. Let $a_{i}^{(m)} = b_{i}^{(m)} + i$, and $a_{i+10^m}^{(m)} =
b_{i+10^m}^{(m)}$, where $1 \le i \le 10^m$. Let $A^{(m)} = \{a_{1}^{(m)},
\dots{} ,a_{10^{m}}^{(m)},a_{10^m+1}^{(m)}, \dots{}
,a_{10^m+10^m}^{(m)}\}$. Define the sets $\displaystyle A = \cup_{m=1}^{\infty}A^{(m)}$ and $B
= \cup_{m=1}^{\infty}B^{(m)}$ and let $A=\{a_1,a_2,\dots \}$, $a_1 < a_2 < \dots $ and $B=\{b_1,b_2,\dots \}$, $b_1 < b_2 < \dots $. It is clear
that $|a_{N} - b_{N}| \ll N$ and $a_{N} = \Theta(N^{3})$ for every $N$, therefore $d_{A,B}(n)\ll n^{1/3}$. In the next step we prove that $A$ may be chosen for a Sidon set. Our strategy is the
following. It is clear from the definition that $A$ is built up from blocks. We use the greedy algorithm to construct the set $A$. Assume that we have already constructed the first few blocks, and we have already chosen some elements to obtain the next block. Suppose that this set satisfies the Sidon property. By the definition of a block, $A^{(m)}$ contains two different type of elements, therefore in each step we have to add two new elements to the set had already been constructed. We have to guarantee the Sidon property. Since we have two new elements, we will need an extra condition to ensure that the sum of the two new elements does not destroy the Sidon property either. More formally, let $B_{0} = \{200, 800\}$ and $A_{0} = \{201, 800\}$ and assume that
we have already chosen blocks $A^{(0)}, A^{(1)}, A^{(2)}, \dots ,A^{(m-1)}$ and integers
$a_{1}^{(m)}, \dots{} ,a_{i-1}^{(m)}, a_{10^m+1}^{(m)}, \dots{} ,a_{10^m+i-1}^{(m)}$, where $1 \le i \le 10^{m}$. For $ 0\le l \le i-1$, define
\[
\mathcal{A}^{(m)}_{l} = \cup_{t=0}^{m-1}A^{(t)} \cup \{a_{1}^{(m)}, \dots{} ,a_{l}^{(m)}, a_{10^m+1}^{(m)}, \dots{} ,a_{10^m+l}^{(m)}\}.
\]
We prove by induction on $l$ that the integers $a_{j}^{(m)}$,  $a_{j+10^m}^{(m)}$ may be selected such that $\mathcal{A}^{(m)}_{l}$ satisfies the Sidon property, that is $a + a^{'} \ne a^{''} + a^{'''}$, where $a, a^{'}, a^{''}, a^{'''} \in \mathcal{A}^{(m)}_{l}$ except for the trivial solutions and the extra condition
\[
a + a^{'} \ne 1000^{m+1} + j,
\]
where $a, a^{'} \in \mathcal{A}^{(m)}_{l}$, and for every $l \le j \le m$.
It is enough to prove that we may add integers $a_{i}^{(m)}$ and $a_{i+10^m}^{(m)}$ to the set $\mathcal{A}^{(m)}_{i-1}$ retaining the Sidon property and the extra condition.
\[
a + a^{'} \ne 1000^{m+1} + j,
\]
where $a, a^{'} \in \mathcal{A}^{(m)}_{i}$, for any $i < j \le 10^{m}$.
In order to
guarantee it we distinguish six cases. Let $a$, $a^{'}$, $a^{''} \in
\mathcal{A}^{(m)}_{i-1}$, and we put $a_{i}^{(m)}$ and $a_{i+10^m}^{(m)}$ to the set $\mathcal{A}^{(m)}_{i-1}$. We have to guarantee that adding new elements  $a_{i}^{(m)}$ and
$a_{i+10^m}^{(m)}$to the set $\mathcal{A}^{(m)}_{i-1}$ does not destroy its Sidon property i.e.,
\begin{equation}
a_{i}^{(m)} + a \ne a^{'} + a^{''},
\end{equation}
and
\begin{equation}
a_{i+10^m}^{(m)} + a \ne a^{'} + a^{''},
\end{equation}
and
\begin{equation}
a_{i}^{(m)} + a_{i+10^m}^{(m)} \ne a + a^{'},
\end{equation}
and
\begin{equation}
a_{i}^{(m)} + a \ne a_{i+10^m}^{(m)} + a^{'}.
\end{equation}
Moreover, by the extra condition we have
\begin{equation}
a_{i}^{(m)} + a \ne 1000^{m+1} + j,
\end{equation}
and
\begin{equation}
a_{10^m+i}^{(m)} + a \ne 1000^{m+1} + j,
\end{equation}
where $i < j \le 10^{m}$.

It is easy to see that the number of elemets of $\mathcal{A}^{(m)}_{i-1}$ is less than $10^{m} + \dots{} + 1$,
thus the number of possibilities to choose triplets $(a,a^{'},a^{''})$ in (6) and (7)
is at most $(10^{m} + \dots{} + 1)^{3} \cdot 2 <
\Big(\frac{10}{9}\cdot 10^{m}\Big)^{3} \cdot 2$.

In the next step we show that inequality (8) holds. It is clear that $a_{i}^{(m)} + a_{i+10^m}^{(m)} = 1000^{m+1} + i$ and the extra condition $a + a^{'} \ne 1000^{m+1} + j$ for any $a, a^{'} \in \mathcal{A}^{(m)}_{i-1}$ and $i - 1 < j \le 10^{m}$. In the special case $j = i$ we have $a + a^{'} \ne 1000^{m+1} + i$, which proves (8).
In equation (9)
\[
a_{i}^{(m)} + a \ne a_{i+10^m}^{(m)} + a^{'} = 1000^{m+1} + i - a_{i}^{(m)} +
a^{'},
\]
implies that
\[
2a_{i}^{(m)} \ne 1000^{m+1} + i + a^{'} - a.
\]
It is clear that the number of possibilities
to choose pairs $(a,a^{'})$ is at most $(1 + \dots{} + 10^{m})^{2} < (\frac{10}{9}10^{m})^{2}$.

It is easy to see that in inequalities (10) and (11) the number of possibilities
to choose pairs $(a,j)$ and is at most $2 \cdot (1 + \dots{} + 10^{m}) \cdot 10^{m} < \frac{20}{9}\cdot (10^{m})^{2}$.

It follows that the number of wrong $a_i^{(m)}$ is at most
\[
2\Big(\frac{10}{9}10^{m}\Big)^{3} + (\frac{10}{9}10^{m})^{2} + \frac{20}{9}\cdot (10^{m})^{2} < 100\cdot 1000^{m}.
\]
This shows that we may choose elements which does not destroy the Sidon
property. $\blacksquare$

\bigskip

\noindent {\bf Proof of Theorem 8.}
It is enough to prove that for every $k\geq 1$, there exists an integer
$n\leq 100^{100^k}$ such that $R_C(n)\geq
k$. To prove this we use the well-known formulas of Vieta: Let us suppose
that $$x^3+y^3=az^3,$$ then we have
$$
(x(x^3+2y^3))^3+(-y(2x^3+y^3))^{3} = a(x^3-y^3)^3z^3.
$$
After repeating we get
$$
\left( x(x^3+2y^3)(x^3(x^3+2y^3)^3-2y^3(2x^3+y^3)^{3})\right)^3
$$
$$
+ \left( y(2x^3+y^3)(2x^3(x^3+2y^3)^3-y^3(2x^3+y^3)^{3})\right)^3
$$
$$
= a(x^3-y^3)^3\left(x^3(x^3+2y^3)^3+y^3(2x^3+y^3)^3\right)^3z^3.
$$
We define the sequences $u_i$, $v_i$ and $w_i$ recursively as
follows. Let $u_1=4^{k-1}$, $v_1=1$ and $w_1=1$ and let
$$
u_{i+1}=u_i(u_i^3+2v_i^3)(u_i^3(u_i^3+2v_i^3)^{3}-2v_i^3(2u_i^3+v_i^3)^{3})
$$
$$
v_{i+1}=v_i(2u_i^3+v_i^3)(2u_i^3(u_i^3+2v_i^3)^{3}-v_i^3(2u_i^3+v_i^3)^{3})
$$
$$
w_{i+1}=(u_i^3-v_i^3)(u_i^3(u_i^3+2v_i^3)^3+v_i^3(2u_i^3+v_i^3)^{3})w_i,
$$
for $i=1,2,\dots ,k-1$. Then we have
$$u_i^3+v_i^3=(64^{k-1}+1)w_i^3$$
and $w_i|w_k$. Let $x_i=\frac{u_i}{w_i}w_k$, $y_i=\frac{v_i}{w_i}w_k$ and
$z_i=w_i$, for $1 \le i \le k$. Hence
$$x_i^3+y_i^3=(64^{k-1}+1)w_k^3.$$
It is enough to show that $0<y_i\leq x_i$, the vectors $(x_i,y_i)$ are
different and $(64^{k-1}+1)w_k^3<100^{100^k}$.
Obviously
\begin{equation}
\frac{x_{i+1}}{y_{i+1}}=\frac{u_{i+1}}{v_{i+1}}=\frac{u_i}{v_i} \cdot
 \frac{u_i^3+2v_i^3}{2u_i^3+v_i^3} \cdot
 \frac{u_i^3(u_i^3+2v_i^3)^{3}-2v_i^3(2u_i^3+v_i^3)^{3}}{2u_i^3(u_i^3+2v_i^3)^{3}-v_i^3(2u_i^3+v_i^3)^{3}}.
\end{equation}
We will show by induction that $u_i\geq 4^{k-i}v_i$ (and therefore
$u_i\geq v_i$) and $u_i,v_i>0$.
This is trivial for $i=1$. Let suppose that $u_i\geq 4^{k-i}v_i\geq 4v_i$ for
some $1 \le i \le k$. Then
we have $$u_i^3(u_i^3+2v_i^3)^{3}-2v_i^3(2u_i^3+v_i^3)^{3} > 0$$ and
therefore by (10) we have
$$\frac{u_{i+1}}{v_{i+1}} \ge \frac{u_i}{v_i} \cdot \frac{1}{2} \cdot
\frac{1}{2},$$ which shows the inductive step.
On the other hand in view of (10) we can see that
$\frac{u_{i+1}}{v_{i+1}}<\frac{u_i}{v_i}$, and therefore we have $k$
different vectors $(x_i,y_i)$.
To finish the proof we have to verify $(64^{k-1}+1)w_k^3<100^{100^k}$.
We know $w_{i+1}\leq 54u_i^{15}w_i$ and therefore
$$w_k\leq 54^{k-1}u_{k-1}^{15}u_{k-2}^{15}\dots u_1^{15}w_1\leq
54^{k-1}u_{k-1}^{15(k-1)}$$
Since $u_{i+1}\leq 81u_i^{16}$, we get
$$u_k\leq 81^{1+16+16^2+\dots +16^{k-2}}\leq 81^{\frac{16^{k-1}}{15}},$$ thus
$$(64^{k-1}+1)w_k\leq (64^{k-1}+1)54^{k-1}81^{(k-1)16^{k-1}}\leq 100^{100^k},$$
which completes the proof. $\blacksquare$

\bigskip

\noindent {\bf Proof of Theorem 9.}
To prove Theorem 9 we use the probabilistic method due to Erd\H{o}s
and R\'enyi. The method is standard therefore we does not give the probabilistic background here (see e.g. the excellent book of Halberstam and Roth [17]).
We denote the probability of an event
by $\mathbb{P}$, and the expectation of a random variable $\zeta$ by
$\mathbb{E}(\zeta)$.
Let $\Omega$ denote the set of the strictly increasing sequences of
positive integers. Theorem 13. in [17], p. 142. shows that one can
obtain a valid probability space ($\Omega$, $X$, $\mathbb{P}$), where
the events
$\mathcal{E}^{(n)} = \{\mathcal{A}$: $\mathcal{A} \in \Omega$, $n \in
\mathcal{A}\}$ are independent for $n = 1, 2, \dots{}$.
We denote the characteristic function of the event $\mathcal{E}^{(n)}$
by $\varrho(\mathcal{A}, n)$:
\[
\varrho(\mathcal{A}, n) =
\left\{
\begin{aligned}
1 \textnormal{, if } n \in \mathcal{A} \\
0 \textnormal{, if } n \notin \mathcal{A}.
\end{aligned} \hspace*{3mm}
\right.
\]
\noindent Furthermore, for some $\mathcal{A} = \{a_1, a_2, \dots{}\} \in
\Omega$ denote by $A(n)$ the number of elements of $\mathcal{A}$ up to
$n$, i.e.,
\[
A(n) = \sum_{\overset{a \in \mathcal{A}}{a \le n}}1.
\]
It is clear that
\[
A(n) = \sum_{j=1}^{n}\varrho(\mathcal{A}, j)
\]
is the sum of Boolean random variables. We need two basic results of
probability theory.
\begin{lem}(Borel-Cantelli)
Let $X_{1}, X_{2}, \dots{}$ be a sequence of events in a probability space. If
\[
\sum_{j=1}^{+\infty}\mathbb{P}(X_{j}) < \infty,
\]
\noindent then with probability 1, at most a finite number of the events
$X_{j}$ can occur.
\end{lem}
\noindent See [17], p. 135. The next tool is the well-known correlation
inequality of Chernoff.
\begin{lem}(Chernoff's inequality)
If $t_i$'s are independent Boolean random
 variables and $X = t_1 + \dots{} + t_n$, then for any $\delta > 0$ we have
$$
\mathbb{P}\big(|X - \mathbb{E}(X)| \ge \delta \mathbb{E}(X)\big) \le
 2e^{-min(\delta^{2}/4, \delta/2)\mathbb{E}(X)}.
$$
\end{lem}
\noindent See in [1]. Define the random sequence $\mathcal{A}$ by
$\mathbb{P}(\{\mathcal{A}$: $\mathcal{A} \in \Omega$, $n \in
\mathcal{A}\}) = \mathbb{P}(n\in \mathcal{A}) =
\frac{1}{3}\frac{1}{n^{2/3}}$ for every positive integer $n$.
It is easy to see that
$$\mathbb{E}(A(x)) = \sum_{n=1}^{x}\frac{1}{3}\frac{1}{n^{2/3}} =
\int_{1}^{x}\frac{1}{3}y^{-2/3}dy + O(1) = x^{1/3} + O(1).
$$ As $A(x)$ is the sum of independent Boolean random
 variables, it follows from Chernoff's inequality with $\delta =
 \frac{3}{x^{1/6}}\sqrt{\log x}$, that
$$
\mathbb{P}\big(|A(x) - \mathbb{E}(A(x))| \ge \frac{3}{x^{1/6}}\sqrt{\log
x}\mathbb{E}(A(x))\big) \le
 2e^{-\frac{9}{4}x^{1/6}\log x \mathbb{E}(A(x))}
$$
$$
< 2e^{-\frac{9}{4}x^{-1/3}(\log x)(x^{1/3}+O(1))} < e^{-2\log x} \leq
\frac{1}{x^{2}},
$$
if $x$ is large enough. By the Borel - Cantelli lemma we have
$$
A(x) = \mathbb{E}(A(x)) + O(x^{1/6}\sqrt{\log x}).
$$
with probability $1$ for every $x \ge 2$. It is clear that
$$
n = A(a_{n}) = {a_{n}}^{1/3} + O(a_{n}^{1/6}\sqrt{\log a_{n}}).
$$
It follows that
$$
a_{n}^{1/3} = n + O(\sqrt{n\log n}),
$$
thus we have
$$
a_{n} = n^{3} + O(n^{5/2}\sqrt{\log n}).
$$
On the other hand we put
\[
r_{2}(\mathcal{A}, n) =  \sum_{1 \le j < n/2}\varrho(\mathcal{A}, j)\varrho(\mathcal{A},
n - j),
\]
which is also a random variable. It is easy to see that
\[
\mathbb{E}(r_{2}(\mathcal{A}, n)) = \frac{1}{9}\sum_{1 \le j <
  n/2}\frac{1}{j^{2/3}(n-j)^{2/3}} \ll \frac{1}{n^{2/3}}\sum_{1 \le j <
  n/2}\frac{1}{j^{2/3}} = \frac{1}{n^{2/3}}\Big(\int_{1}^{n/2}j^{-2/3}dj +
O(1)\Big) \ll n^{-1/3}.
\]
Let $E_{i}$ be the event
\[
E_{i} = \{i \in \mathcal{A}, n - i \in \mathcal{A}\}.
\]
It is clear that the events $E_{i}$'s are mutually independent. Thus we have
\[
\mathbb{P}(r_{2}(\mathcal{A}, n) > 3) \le  \sum_{1 \le i_{1} < i_{2} <
i_{3} < i_{4} < n/2}\mathbb{P}(E_{i_{1}} \cap \dots{} \cap
E_{i_{4}}) = \sum_{1 \le i_{1} < i_{2} < i_{3} < i_{4} <
    n/2}\mathbb{P}(E_{i_{1}}) \dots \mathbb{P}(E_{i_{4}})
\]
\[
\le \Big(
\frac{1}{9}\sum_{1 \le j < n/2}\frac{1}{j^{2/3}(n-j)^{2/3}}\Big)^{4}
= \mathbb{E}(r_{2}(\mathcal{A}, n))^{4} \ll n^{-4/3}.
\]
It follows from the Borel - Cantelli lemma that with probability $1$,
for $n > n_{0}$, $r_{2}(\mathcal{A}, n) \le 3$. This
implies that $R_{\mathcal{A}}(n)$ is also bounded and so does
$s_{\mathcal{A}}(n)$.  $\blacksquare$


\begin{thebibliography}{99}
\bibitem{als} N. Alon - J. Spencer, \textit{The Probabilistic method},
	Wiley Interscience (2000).
\bibitem{cil} P. Borwein, S. Choi, F. Chu, \textit{An old conjecture of
  Erd\H{o}s - Tur\'an on additive bases}, Math. Comp. \textbf{75} (2006),
  475-484.
\bibitem{des} Y-G. Chen, \textit{The analogue of Erd\H{o}s - Tur\'an
  conjecture in $\mathbb{Z}_{m}$}, J. Number Theory, \textbf{128} (2008),
  2573-2581.
\bibitem{gao} Y-G. Chen, \textit{On the Erd\H{o}s - Tur\'an
  conjecture}, C. R. Math. Acad. Sci. Paris \textbf{350} (2012),
  933-935.
\bibitem{ert} G. A. Dirac, \textit{Note on a problem in additive number
  theory}, J. London Math. Soc., \textbf{26} (1951), 312-313.
\bibitem{est} M. Dowd, \textit{Questions related to the Erd\H{o}s - Tur\'an
  conjecture }, SIAM J. Discrete Math., \textbf{1} (1988), 142-150.
\bibitem{hal} P. Erd\H{o}s, P. Tur\'an, \textit{On a problem of Sidon in
  additive number theory, and some related problems}, J. London Math. Soc.
  \textbf{16} (1941), 212-215.
\bibitem{kim} P. Erd\H{o}s, \textit{Problems and results in
  additive number theory}, in: Colloque sur la Th\'eorie des Nombres,
  Bruxelles, 1955, George Thone, Li\'ege; Masson and Cie, Paris, 1956,
  pp. 127-137.
\bibitem{elk} P. Erd\H{o}s, W. H. J. Fuchs, \textit{On a problem of additive
  number theory}, J. London Math. Soc., \textbf{31} (1956), 67-73.
\bibitem{ghl} G. Grekos, L. Haddad, C. Helou, J. Pihko,
    \textit{On the general Erd\H{o}s - Tur\'an conjecture},
	International Journal of Combinatorics, (2014) Art. ID 826141,
	11 pp.
\bibitem{ghh} G. Grekos, L. Haddad, C. Helou, J. Pihko,
    \textit{On the Erd\H{o}s - Tur\'an conjecture}, J. Number Theory
    \textbf{102} (2003), 339-352.
\bibitem{ghp} G. Grekos, L. Haddad, C. Helou, J. Pihko,
    \textit{The class of Erd\H{o}s - Tur\'an sets}, Acta Arith., \textbf{117}
    (2005), 81-105.
\bibitem{heh} G. Grekos, L. Haddad, C. Helou, J. Pihko,
    \textit{Variations on a theme of Cassels for additive bases},
    Intl. J. Number Theory \textbf{2} (2006), 249-265.
 \bibitem{hel} L. Haddad, C. Helou,
    \textit{Bases in some additive groups and the Erd\H{o}s - Tur\'an
      conjecture}, J. Combin. Theory Ser. A
    \textbf{108} (2003), 339-352.
\bibitem{hec} L. Haddad, C. Helou,
    \textit{Additive bases representations in groups}, Integers
    \textbf{8} (2008), A5.
\bibitem{het} L. Haddad, C. Helou,
    \textit{Representations of integers by near quadratic sequences}, Journal
    of Integer sequences
    \textbf{15} (2012), 12.8.8.
\bibitem{hal} H. Halberstam, K. F. Roth, \textit{Sequences}, Springer -
  Verlag, New York, 1983.
\bibitem{kov} S. F. Konyagin and V. F. Lev, \textit{The Erd\H{o}s - Tur\'an
  problem in infinite groups}, in: Additive Number Theory, Springer \textbf{20} (2010), 195 - 202.
\bibitem{nab} M. B. Nathanson, \textit{Unique representation bases for the
  integers}, Acta Arith., \textbf{108} (2003), 1 - 8.
 \bibitem{nar} M. B. Nathanson, \textit{Representation functions of additive
  bases for abelian semigroups}, Int. J. Math. Sci. (2004), 1589  - 1597.
\bibitem{nat} M. B. Nathanson, \textit{Generalized additive bases, K\H{o}nig's
lemma and the Erd\H{o}s - Tur\'an conjecture}, \textbf{106} (2004), 70 - 78.
\bibitem{son} J.Nesetril and O. Serra, \textit{On a conjecture of Erd\H{o}s and
Tur\'an for additive basis}, in Proceedings of the ``Segundas Jornadas de
  Theoria de N\'umeros'', Bibl. Rev. Mat. Iberoamericana, (2008), 209 - 220.
\bibitem{niz} I. Niven, H. S. Zuckerman, H. L. Montgomery, \textit{An Introduction to The Theory of Numbers}, 5th ed., Wiley, 1991.
\bibitem{tav} V. Pus, \textit{On multiplicative bases in abelian groups},
  Czech. Math. J. \textbf{41} (1991), 282 - 287.
\bibitem{val} Cs. S\'andor, \textit{A note on a conjecture of Erd\H{o}s and
Tur\'an}, Integers, \textbf{8} (2008), A30.
\bibitem{tan} T. Min, \textit{On the Erd\H{o}s - Tur\'an conjecture},
  J. of Number Theory, \textbf{150} (2015), 74 - 80.
\bibitem{teb} G. Tenenbaum, \textit{Introduction to analytic and probabilistic number theory}, Cambridge University Press; 1st ed. 1995.
\bibitem{tav} Y. Quan-Hui, \textit{A generalization of Chen's theorem on
	the Erd\H{o}s-Tur\'an conjecture}, Int. J. Number Theory
	\textbf{9} (2013), 1683 - 1686.


\end{thebibliography}
\end{document}